\documentclass[conference]{IEEEtran}

\usepackage{tabularx}
\usepackage{hyperref}
\usepackage[english]{babel}
\usepackage{amsmath} 
\usepackage{amssymb}  
\usepackage{graphicx}
\usepackage{balance}
\usepackage[colorinlistoftodos]{todonotes}
\newtheorem{remark}{Remark}[section]

\ifCLASSINFOpdf
\else
\fi
\let\svthefootnote\thefootnote

\usepackage{amsmath}
\begin{document}
%
\title{Model-free Load Control for High Penetration of \\ Solar Photovoltaic Generation }

%
%
%


%
\author{\IEEEauthorblockN{Ouassim Bara\IEEEauthorrefmark{1},
Mohammad Olama\IEEEauthorrefmark{4}, Seddik Djouadi\IEEEauthorrefmark{1}, Teja Kuruganti\IEEEauthorrefmark{4}, Michel Fliess\IEEEauthorrefmark{2}\IEEEauthorrefmark{5} and
C\'{e}dric Join\IEEEauthorrefmark{3}\IEEEauthorrefmark{5}\IEEEauthorrefmark{6}
 }
\IEEEauthorblockA{\IEEEauthorrefmark{1}Department of Electrical Engineering \& Computer Science, University of Tennessee,
Knoxville, TN 37996, USA. \\ Email: \{obara, mdjouadi\}@utk.edu}
\IEEEauthorblockA{\IEEEauthorrefmark{4}Oak Ridge National Laboratory,	Oak Ridge, TN 37831, USA. Email:\{olamahussemm, kurugantipv\}@ornl.gov}
\IEEEauthorblockA{\IEEEauthorrefmark{2}LIX (CNRS, UMR 7161), \'Ecole polytechnique, 91128 Palaiseau, France. Email: Michel.Fliess@polytechnique.edu}
\IEEEauthorblockA{\IEEEauthorrefmark{3}CRAN (CNRS, UMR 7039), Universit\'{e} de Lorraine, BP 239, 54506 Vand{\oe}uvre-l\`{e}s-Nancy, France. \\ Email: cedric.join@univ-lorraine.fr}

\IEEEauthorblockA{\IEEEauthorrefmark{5}AL.I.E.N., 7 rue Maurice Barr\`{e}s, 54330 V\'{e}zelise, France. Email: \{michel.fliess, cedric.join\}@alien-sas.com}
\IEEEauthorblockA{\IEEEauthorrefmark{6}Projet Non-A, INRIA Lille -- Nord-Europe, France}}


\maketitle

\begin{abstract}
This paper presents a new model-free control (MFC) mechanism that enables the local distribution level circuit consumption of the photovoltaic (PV) generation by local building loads, in particular, distributed heating, ventilation and air conditioning (HVAC) units. The local consumption of PV generation will help minimize the impact of PV generation on the distribution grid, reduce the required battery storage capacity for PV penetration, and increase solar PV generation penetration levels. The proposed MFC approach with its corresponding intelligent controllers does not require any precise model for buildings, where a reliable modeling is a demanding task. Even when assuming the availability of a good model, the various building architectures would compromise the performance objectives of any model-based control strategy. The objective is to consume most of the PV generation locally while maintaining occupants comfort and physical constraints of HVAC units. That is, by enabling proper scheduling of responsive loads temporally and spatially to minimize the difference between demand and PV production, it would be possible to reduce voltage variations and two-way power flow. Computer simulations show promising results where a significant proportion of the PV generation can be consumed by building HVAC units with the help of intelligent control.

\end{abstract}

\begin{IEEEkeywords} 
HVAC, distributed energy resources, building load control, solar variability, model-free control, intelligent proportional controllers. 
\end{IEEEkeywords}


%
\IEEEpeerreviewmaketitle

\section{Introduction}
A large introduction of solar photovoltaic (PV) and other renewable resources may significantly alter the stability and balance of the power grid, through voltage fluctuations and the ability to maintain the required network frequency \cite{kim2013}. A viable solution to this problem, knowing the widespread use of  HVAC systems in residential and commercial buildings, is to use the idea of load shaping. By introducing responsive loads (buildings), the control of HVAC units will act as a counter to undesirable fluctuations of power generated by the building sized PV arrays, while not sacrificing some level of temperature comfort. 

The motivation of not relying on model-based control approaches comes from the daunting task of devising a reliable building model, either via physical laws or via blackbox identification, and also,  in the fact that it is difficult to account for all the uncertainties and unknown disturbances, especially with respect to strong weather disturbances and changes in occupancy pattern (see, \textit{e.g.}, \cite{nz,mont,teodorescu20111,raoufet,zhou} and the references therein). A good model often requires the installation of too many sensors necessary for the control algorithm decision making part. Not to mention the financial cost that this will incur. \let\thefootnote\relax\footnote{This manuscript has been authored by UT-Battelle, LLC under Contract No. DE-AC05-00OR22725 with the U.S. Department of Energy.
}\addtocounter{footnote}{-1}\let\thefootnote\svthefootnote

We are therefore proposing to use a new model-free control (MFC) approach with its corresponding intelligent controllers \cite{csm}, where the need of any precise modeling disappears. This setting 
\begin{itemize}
\item is \emph{data-driven}, \textit{i.e.}, only the input and output data are used. Mathematical physical laws and the associated differential equations are ignored.
\item is not limited to finite-dimensional linear systems. It has been successfully applied in much more involved situations with nonlinearities and where modeling via partial differential equations has been proposed.\footnote{Employing PDEs is quite common in the air conditioning of buildings. See, \textit{e.g.}, \cite{pde} and the references therein.}
\item is easy to implement and quite robust with respect to internal and external disturbances, even when compared to the most popular Proportional Integral Derivative controller (see, \textit{e.g.}, \cite{astrom});
\item has already been successfully employed all over the world and gave rise to some patents: see, \textit{e.g.}, the references in \cite{csm}, \cite{alinea} and the references therein, and \cite{eccDenmark,buc,portugal,maa,panc,its,kr}. For obvious reasons let us emphasize here  the heating and the humidification of an agricultural greenhouse \cite{toulon}, and the heating of a single building \cite{bldg}. 
\end{itemize}

Although the idea of using HVAC unit systems to provide ancillary services to the grid is not new (see \cite{hao2014}, \cite{balanda2014}), the  main contribution of this paper is on designing a MFC framework in order to consume most of the generated PV power locally by  building loads, while simultaneously ensuring a certain temperature comfort inside the buildings. How to solve  the actual problem turns out to be  intuitive, in the sense that one  can consider the total energy provided by the PV profile and share it amongst the available building load. A straightforward solution was then to consider at each instant of time the available energy coming from the profile as a constraint for each HVAC system. 
Note that the difficulty in trying to implement such procedures is the necessity for the control to follow two reference trajectories simultaneously (see also \cite{siaap}). One corresponding to the interior temperature of the buildings (set in advance within a certain comfort zone), and the other, corresponding to the generated PV profile.  It is clearly a hard requirement to satisfy if one wants to avoid  optimization problem, where it is easier to introduce those constraints. The price to pay for such optimization will result in  higher computational power, even if sub-optimality mathematical techniques are introduced (see, \textit{e.g.}, \cite{pde,tn,zheng,swarm} and the references therein). A real-time implementation would thus become problematic. 

The remainder of the paper is organized as follows.  Section \ref{sec:mfc} provides a short introduction to model-free control control and to the corresponding intelligent controllers.\footnote{See \cite{csm} for full details.} 
 A simplified physical model is presented in Section \ref{sec:formulation} together with the formulation of our approach and the simulation results. Section \ref{conclusion} summarizes the paper and presents the conclusions.

\section{Model-free control and intelligent controllers} \label{sec:mfc}

Without any loss of generality we restrict ourselves to Single-Input Single-Output systems. Instead of trying to write down a complex differential equation, introduce the \emph{ultra-local} model 
\begin{equation}
{\dot{y} = F + \alpha u} \label{1}
\end{equation}
where
\begin{itemize}
	\item $u$ and $y$ are the input (control) and output variables,
	\item the derivation order of $y$ is $1$, like in most concrete situations,
	\item $\alpha \in \mathbb{R}$ is chosen by the practitioner such that 
	$\alpha u$ and $\dot{y}$ are of the same magnitude.
\end{itemize}
The following explanations on $F$ might be useful: 
\begin{itemize}
\item $F$ subsumes the knowledge of any model uncertainties and disturbances,
	\item $F$ is estimated via the measures of $u$ and $y$.
	\end{itemize}

\subsection{Intelligent controllers}
The loop is closed by an \emph{intelligent proportional controller}, or \emph{iP},
\begin{equation}\label{ip}
u = - \frac{\hat{F} - \dot{y}^\ast + K_P e}{\alpha}
\end{equation}
where
\begin{itemize}
	\item $y^\star$ is the reference trajectory,
	\item $e = y - y^\star$ is the tracking error,
	\item $K_P$ is the usual tuning gain.
\end{itemize}
Combining equations \eqref{1} and \eqref{ip} yields:
$$
	\dot{e} + K_P e = 0
$$
where $F$ does not appear anymore. Local exponential stability is ensured if $Kp>0$: 
\begin{itemize}
\item The gain $K_P$ is thus easily tuned. 
\item Robustness with respect to different types of disturbances and model uncertainties is achieved.
\end{itemize}
\subsection{Estimation of $F$}\label{F}
$F$ is estimated in real-time  according to recent algebraic identification techniques \cite{sira1,sira2,sira3}. 

\subsubsection{First approach}
The term $F$ in Equation \eqref{1} may be assumed to be ``well'' approximated by a piecewise constant function $\hat{F} $. Rewrite then Equation \eqref{1}  in the operational domain (see, \textit{e.g.}, \cite{yosida}): 

\begin{align}
	sY = \frac{\Phi}{s}+\alpha U +y(0)
\end{align}

\noindent where $\Phi$ is a constant. We get rid of the initial condition $y(0)$ by multiplying both sides on the left by $\frac{d}{ds}$:
\begin{align}
Y + s\frac{dY}{ds}=-\frac{\Phi}{s^2}+\alpha \frac{dU}{ds}
\end{align}
Noise attenuation is achieved by multiplying both sides on the left by $s^{-2}$. It yields in the time domain the real-time estimate, thanks to the equivalence between $\frac{d}{ds}$ and the multiplication by $-t$,
\begin{equation}\label{integral1}
{\small \hat{F}(t)  =-\frac{6}{\tau^3}\int_{t-\tau}^t \left\lbrack (\tau -2\sigma)y(\sigma)+\alpha\sigma(\tau -\sigma)u(\sigma) \right\rbrack d\sigma }
\end{equation}

\subsubsection{Second approach}\label{2e}
Close the loop with the iP \eqref{ip}:
\begin{equation}\label{integral2}
\hat{F}(t) = \frac{1}{\tau}\left[\int_{t - \tau}^{t}\left(\dot{y}^{\star}-\alpha u
- K_P e \right) d\sigma \right] 
\end{equation}
\begin{remark}
	Note the following facts: 
	\begin{itemize}
		\item Integrals \eqref{integral1} and \eqref{integral2} are low pass filters.
		\item $\tau > 0$ may be chosen quite small.
		\item Integrals \eqref{integral1} and \eqref{integral2} may of course be replaced in practice by classic digital filters.
	\end{itemize}
\end{remark}

\section{Problem formulation and computer simulations} \label{sec:formulation}
\subsection{A simple mathematical model}\label{sec:mathmodel}
Computer simulations require obviously some mathematical modeling. We have selected a simple set of linear differential equations, which are time-invariant, \textit{i.e.}, with constant coefficients. It was derived in \cite{gwerder2005} and successfully used in \cite{tn,oldewurtel2010,aut}. The dynamic of the interior temperature, interior wall surface temperature, and exterior wall core temperature are given by
\begin{align}
\dot{T}_1&=\frac{1}{C_1}\big[(K_1+K_2)(T_2-T_1)+K_5(T_3-T_1) +u_c \nonumber \\
&+\delta_2+\delta_3 \big] \nonumber \\
\dot{T}_2 &= \frac{1}{C_2}\big[ (K_1+K_2)(T_1-T_2)+\delta_2\big] \nonumber \\
\dot{T}_3 &= \frac{1}{C_3}\big[K_5(T_1-T_3)+K_4(\delta_1-T_3)\big] 
\end{align}
\begin{align}\label{stsp}
\dot{x}=A_cx+B_cu+C_cw 
\end{align}
where the state vector $x$, the control $u$ and the disturbance $w$ are given respectively as
\begin{align}
x=\begin{bmatrix} T_1\\T_2\\T_3 \end{bmatrix}, \quad u=u_c,\quad w=\begin{bmatrix}\delta_1\\ \delta_2\\ \delta_3	\end{bmatrix} \nonumber
\end{align}
$T_1$: room air temperature [$^\circ C$]\hspace{7mm} $ 22\leq T_1\leq 24$,\\
$T_2$: interior-wall surface temperature [$^\circ C$], \\
$T_3$: exterior-wall core temperature [$^\circ C$], \\
$u_c$: cooling power ($\leq 0$) [$kW$], \\
$\delta_1$: outside air temperature [$^\circ C$],\\
$\delta_2$: solar radiation [$kW/m^2$], \\
$\delta_3$: internal heat sources [$kW$],\\
$C_1 = 9.356 \times 10^5$ kJ/C, $C_2 = 2.970 \times 10^6$ kJ/C, $C_3 = 6.695 \times 10^5$  kJ/C, $K_1 = 16.48$ kW/C, $K_2 = 108.5$ kW/C, $K_3 = 5$, $K_4 = 30.5$ kW/C and $K_5 = 23.04$ kW/C.

The above equations are only valid for a summer day. The main control variable is therefore cooling power $u_c$.
One of our control objectives is to regulate the interior temperature of several buildings around 23 $^\circ C$.  This is the reference trajectory.
\begin{remark}
A poor knowledge of the coefficients in Equation \eqref{stsp} is unavoidable in practice. Lack of space prevents us of confirming via simulations the robustness of our MFC setting in such a situation.
\end{remark}

\subsection{Temperature regulation and PV following}
Our aim is to regulate the interior temperature of multiple buildings while at the same time assuring the whole consumed energy tracks as closely as possible the PV energy profile. 
 It yields the following control objectives:
\begin{itemize}
	\item Maintain the interior temperature under an acceptable comfort zone, \textit{i.e.}, between $22^\circ C$ and $24^\circ C$.
	\item Follow the generated PV profile as closely as possible.
\end{itemize}  
\begin{remark}
The PV profile is the only energy that the HVAC units are supposed to get. The buildings are also connected to the grid. It means that the total energy consumed is allowed to fluctuate around the generated PV signal.
\end{remark}

Instead of considering the PV profile as a constraint, as it has been mentioned in the introduction, \textit{i.e.}, 
\begin{align}
\sum_{i=1}^{N_b}u_i(t) \leq PV(t)
\end{align}
where $N_b$ is the number of buildings, a better alternative is to assume a band around the PV profile, such that the available energy  is given according to
\begin{align} \nonumber
E_d(t):\begin{cases} E_d(t) = 0 \quad \text{if} \quad PV(t) = 0 \\ PV(t)-\epsilon \leq E_d(t) \leq PV(t)+\epsilon \quad \text{if} \quad PV(t)> 0\end{cases}
\end{align}
where 
\begin{itemize}
\item $\epsilon > 0$ is a constant, which is chosen according to the grid capabilities and the comfort resulting from the interior temperature;
\item it is assumed that, when $E_d(t)=0$, all the buildings are connected to the grid and no PV tracking is performed. 
\item The constraint $0\leq u_c \leq 3$ on the cooling power must be satisfied when there is no PV energy.
\end{itemize}
Therefore the total energy consumption is 
\begin{align} \label{eq:ineq1}
PV(t) -\epsilon\leq \sum_{i=1}^{N_b}u_i(t) \leq PV(t) +\epsilon 
\end{align}
Assume for simplicity's sake that all the buildings are identical. Then Inequality \eqref{eq:ineq1} yields
\begin{align} \label{eq:ineq2}
\frac{PV(t)}{N_b}-\frac{\epsilon}{N_b}\leq u_i (t) \leq \frac{PV(t)}{N_b}+\frac{\epsilon}{N_b}, \quad  i=1,2,..N_b
\end{align}


\subsection{Computer simulations}
A $10$ minute sampling is used.\footnote{The time lapse depends heavily, of course, on the nature of the plant we want to regulate.} Consider first 13 buildings. 
The simulations are presented first when only the interior temperatures $T{_1}_i$, $i=1,\dots N_b$, is regulated. Afterwards a second objective is added, \textit{i.e.},  the tracking of the PV generated signal, \textit{i.e.}, Inequality \eqref{eq:ineq1} must be satisfied.

\subsubsection{Control design}
Since the interior temperatures $T{_1}_i$ are available for measurements, Equation \eqref{1} yields the ultra-local model
\begin{align}
	\dot{T}{_1}_i = F_i + \alpha_i u{_c}_i 
\end{align} 
From Equation (\ref{ip}) we deduce the following intelligent proportional controller
\begin{align}
	u{_c}_i = -\frac{\hat{F}_i - \dot{T_1}^*_i +K{_P}_i e{_p}_i}{\alpha_i} 
\end{align}
where 
\begin{itemize}
\item $\alpha_i = 5$, $K{_P}_i  = 2$,
\item ${e_p}_i = T{_1}_i -{T_1}^*_i$ is the tracking error,
\item the reference trajectory ${T_1}^*_i$  should be as close as possible to the reference temperature $23^\circ C$,
\item $\hat{F}_i$ is an estimate of $F$, which is obtained according to Section \ref{F},
\item for a ``good'' PV tracking each control $u{_c}_i$ must satisfy Inequality (\ref{eq:ineq2}).
\end{itemize}

\subsubsection{Model-free control with no PV tracking}
\begin{figure}[!t]
\hspace{-2mm} \includegraphics[scale=0.46]{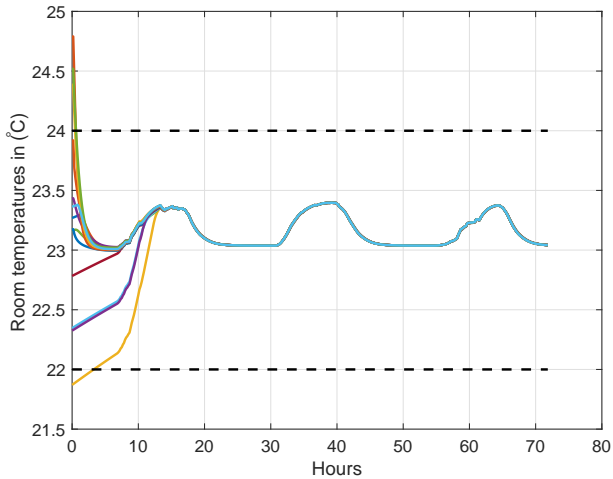}
	\caption{Interior temperature variations in $^\circ C$ for  $N_b=13$} 
	\label{Fig:tempnotracking}
\end{figure}

\begin{figure}[t]
	\hspace{-3mm} \includegraphics[scale=0.46]{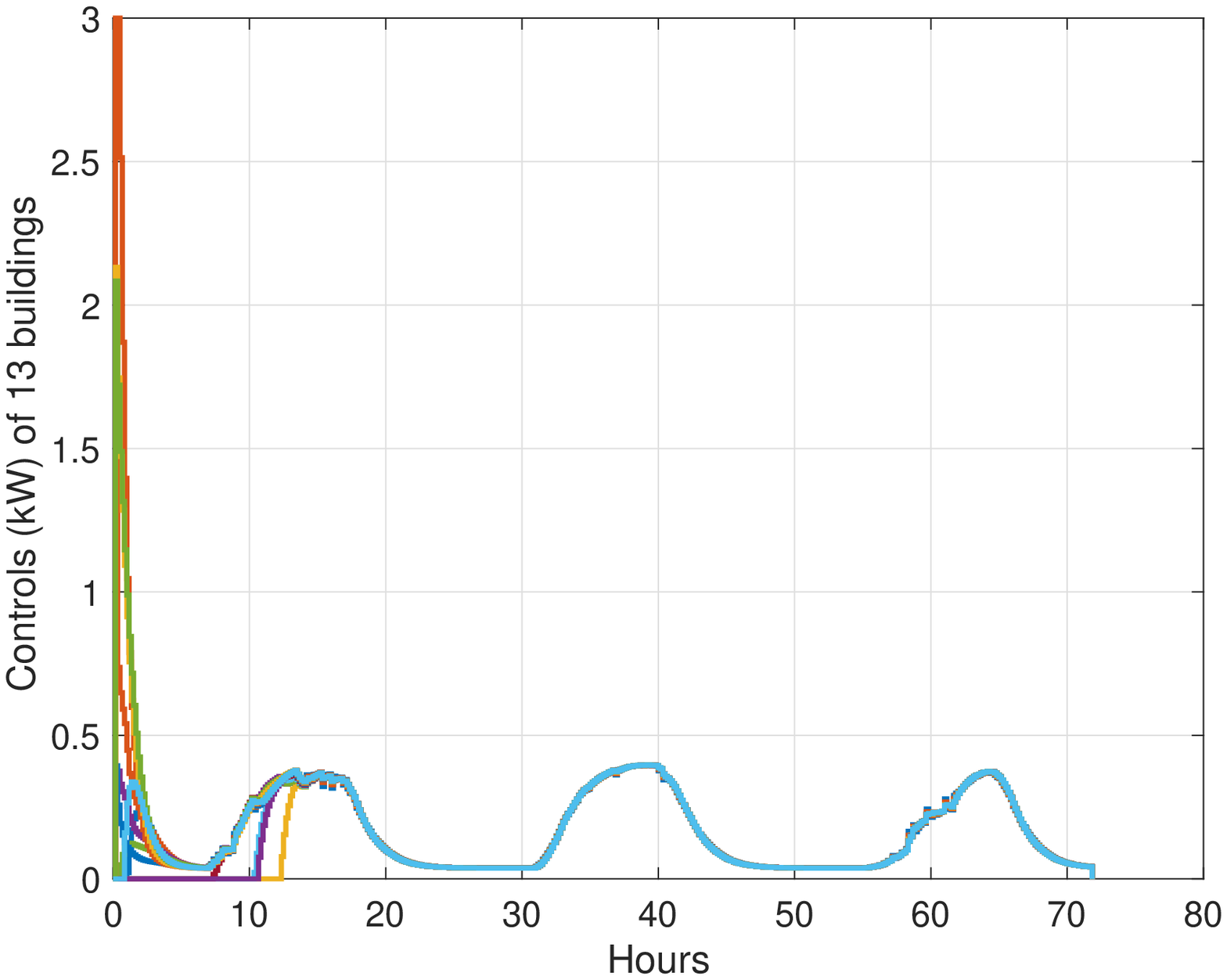} 
	\caption{Time evolution of the $13$ control inputs }  \label{Fig:controlnotracking}
\end{figure}

\begin{figure}[ht]
	\hspace{-2mm} \includegraphics[scale=0.46]{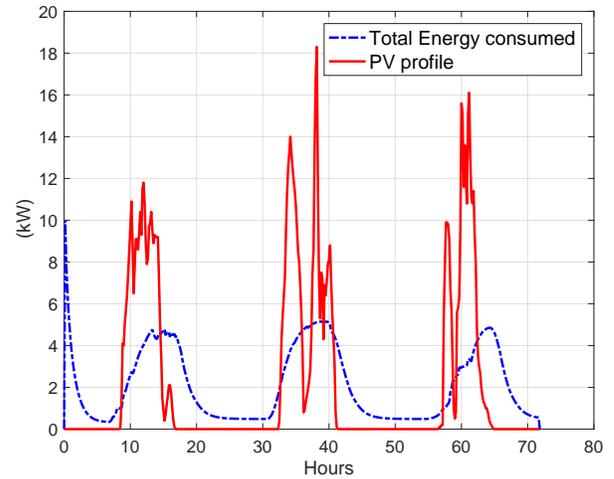} 
	\caption{Total energy $\sum_{i=1}^{N_b}{u_c}_i(t) $  and the PV profile} \label{Fig:pvnotracking}
\end{figure}

Figure \ref{Fig:tempnotracking} displays the evolution of the interior temperature during three days. The temperature is perfectly regulated. MFC is very efficient in rejecting outside  disturbances: one can see their effect in the little bumps appearing essentially from the $8^{th}$ to the $20^{th}$ hour of the day. Note also that the initial temperatures in each building are different. The control inputs in Figure \ref{Fig:controlnotracking} satisfy the required constraints. The initial large values correspond to the buildings with the highest initial temperature. It is important to notice in Figure \ref{Fig:pvnotracking} that the total energy used by the HVAC units does not take into account the variations of the PV profile, hence the large deviation of the  total energy signal represented in blue dashed lines from the actual PV given in red solid line.


\subsubsection{Simulation with PV profile}
\begin{figure}
	\hspace{-2mm} \includegraphics[scale=0.45]{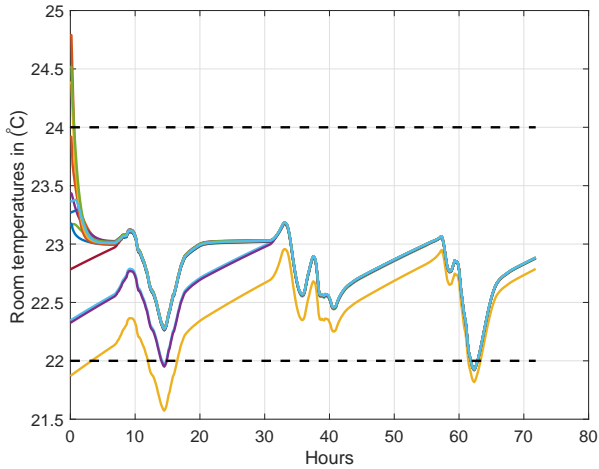}
	\caption{Interior temperature variations in $^\circ C$ for  $N_b=13$} \label{Fig:temptracking}
\end{figure}

\begin{figure}[t]
	\hspace{-4mm} \includegraphics[scale=0.45]{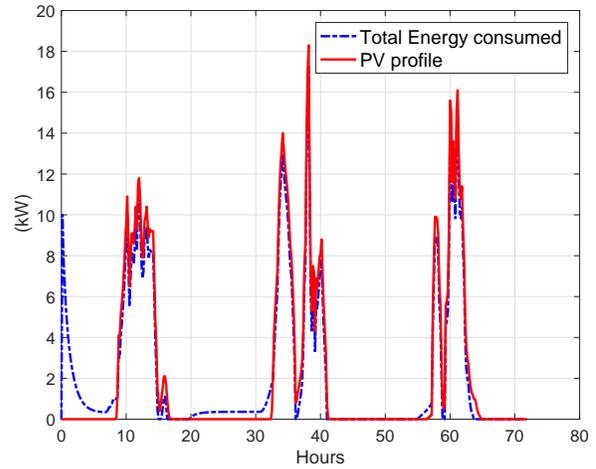} 
	\caption{Total energy $\sum_{i=1}^{N_b}{u_c}_i(t) $  and the PV profile ($N_b = 13$)}   \label{Fig:pvtracking}
\end{figure}

As soon as the photovoltaic source starts collecting energy, Inequality (\ref{eq:ineq2}) is taken into account.  Figure \ref{Fig:temptracking} displays the time evolution of the interior temperature for 13 buildings while the PV tracking is taken into account. Note that most of the temperatures are within their prescribed comfort zone. The slight violation of the lower constraints for one or two buildings is due to the tight restriction on the PV tracking as shown by Figure \ref{Fig:pvtracking}:  we select $\epsilon = 1$, \textit{i.e.}, a rather narrow tolerance margin. Temperatures where slight overshoots appear are those that are starting from the lowest initial temperatures. By  imposing more energy consumption for the buildings than necessary for regulating its indoor temperature, one would expect the HVAC units to run during more time or with higher intensity. It would yield large variations. This behavior is rather normal. It can be significantly improved by selecting an appropriate number of buildings. It is important to understand the difficulty of the actual compromise between the requirement of imposing a certain level of comfort inside the buildings, and that of closely following the generated PV signal. Those two requirements are very often contradictory. Our viewpoint helps in satisfying both control objectives without the need of any complex optimization procedure, especially if one has to deal with a  nonlinear model. As an initial study we decided to focus on a limited number of HVAC units to assess our results. However scaling the problem to larger building units is straightforward. It does not necessitate any change in the problem formulation. More importantly it will not require any additional significant computational power. 

As already mentioned, it is possible to improve the results in Figure \ref{Fig:temptracking}. Add for that purpose one building, \textit{i.e.}, $N_b=14$. Figure \ref{Fig:temptracking14} shows a clear-cut improvement with respect to the interior temperature. Figures 8 and 9  do not display on the other hand any significant change with respect to the case of 13 buildings. One needs to be cautious when dealing with an unfeasible scenario where the maximum energy consumed by all buildings is very small when compared to the energy coming from the PV source, or when the PV energy is not sufficient to account for all the HVAC units. Thus, a better temperature regulation emphasizes that we are, in fact, in the scenario where more energy was available from the PV  for the 13 buildings. Adding one more helped in distributing the extra energy to the new building, putting less weight on the remaining ones to follow the PV profile.
 
 The above results are quite good. Remember that enabling the buildings to absorb the PV sources will help the latter to minimize their negative effect  on the distribution grid and consequently to provide responsive loads that would ultimately help reducing voltage variations and two-way power flow.
 

\section{Conclusion} \label{conclusion}
In order to offer a better integration of the photovoltaic energy sources into the grid, a model-free control approach has been proposed. Computer simulations are showing that the responsive loads (HVAC units)  may be regulated in such a way that the tracking error between demand and PV generation is minimized. This fact will help reducing voltage variations and two-way power flow. The two main control objectives mentioned in the introduction are thus satisfied:
\begin{itemize}
\item the interior temperatures inside each unit remains in the band $\pm1^\circ C$  around the reference temperature,
\item the tracking error of the PV profile is less than $1$kW. 
\end{itemize}
A forthcoming publication will extend this work to Boolean controls,  \textit{i.e.}, on/off controls. This step should lead to concrete implementations. Let us emphasize finally that the hardware implementation of our control strategy is cheap and easy \cite{nice}.

\section*{ACKNOWLEDGMENT}
This material is based upon work supported by the U.S. Department of Energy, Office of Energy Efficiency \& Renewable Energy, the SunShot National Laboratory Multiyear Partnership (SuNLaMP) program under Contract No. DE-AC05-00OR22725.


\begin{figure}
	\hspace{-2mm} \includegraphics[scale=0.45]{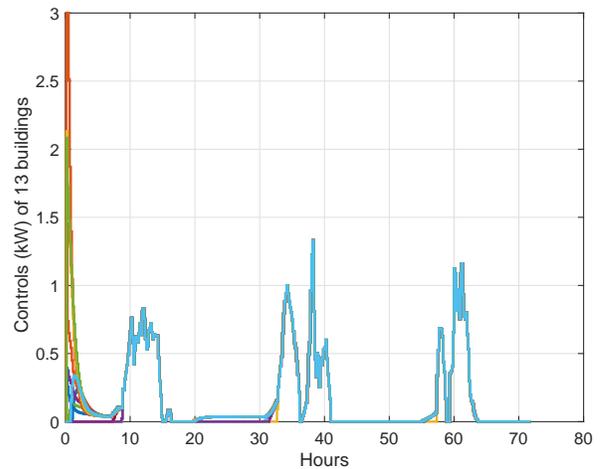} 
	\caption{Time evolution of the $13$ control inputs}   \label{Fig:controltracking}
\end{figure}


\balance

\begin{figure}[t!]
	\hspace{-4mm} \includegraphics[scale=0.45]{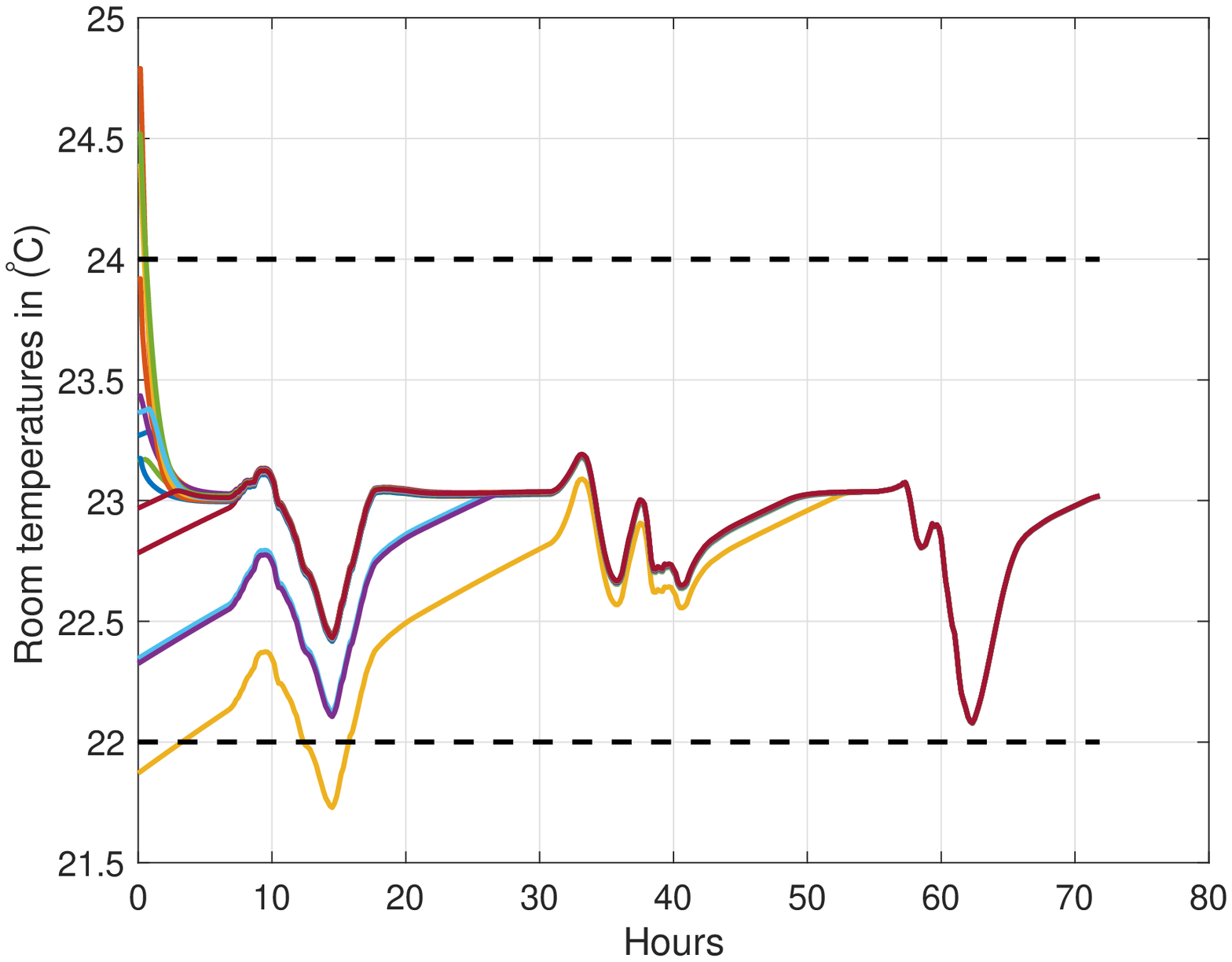}
	\caption{Interior temperature variations in $^\circ C$ for  $N_b=14$} \label{Fig:temptracking14}
\end{figure}

\begin{figure}[t!]
	\hspace{-4mm} \includegraphics[scale=0.45]{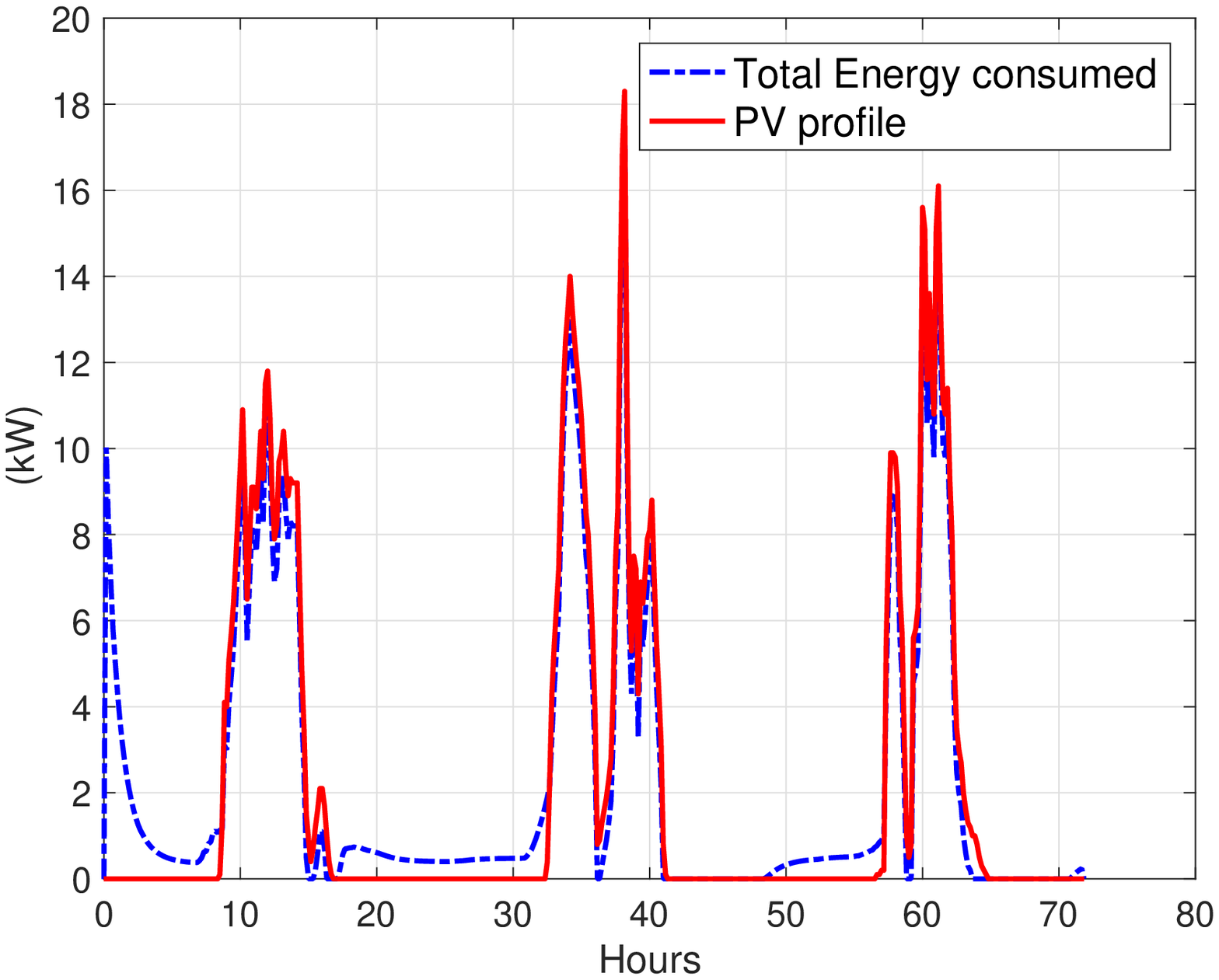} 
	\caption{Total energy $\sum_{i=1}^{N_b}{u_c}_i(t) $  and the PV profile ($N_b=14$)}\label{Fig:pv14}   
\end{figure}

\begin{figure}[t!]
	\hspace{-4mm} \includegraphics[scale=0.45]{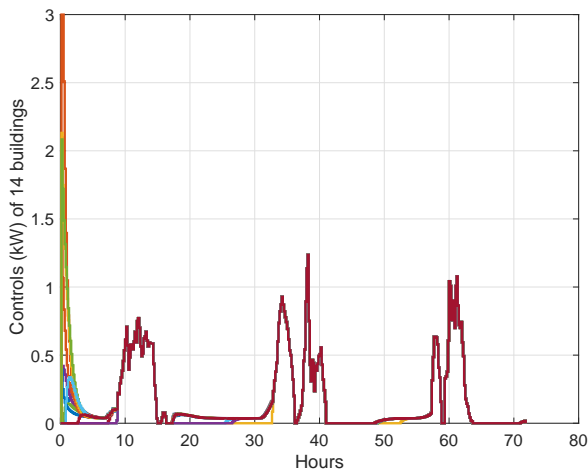} 
	\caption{Time evolution of the $14$ control inputs}   \label{Fig:controltracking14}
\end{figure}


\begin{thebibliography}{9}
	
	\bibitem{kim2013}
	Y.-S. Kim, I.-Y. Chung, S.-I. Moon. 
	\newblock An analysis of variable-speed wind turbine power-control methods with fluctuating wind speed.
	\newblock \emph{Energies}, 6, 3323-3338, 2013.

\bibitem{nz}
D. Jovcic, N. Pahalawaththa, M. Zavahir. Analytical modelling of HVDC-HVAC systems. \emph{IEEE Trans. Power Deliv.}, 14, 506-511, 1999.

\bibitem{mont}
R. Montgomery, R. McDowall. \emph{Fundamentals of HVAC Control Systems}. Elsevier, 2009.

\bibitem{teodorescu20111}
R. Teodorescu, M. Liserre, P. Rodriguez.
\newblock \emph{Grid Converters for Photovoltaic and Wind Power Systems}.
\newblock Wiley, 2011.

    \bibitem{raoufet}
    M.E. Raoufat, A. Khayatian. A. Mojallal.
    \newblock Performance recovery of voltage source converters with application to grid-connected fuel cell DGs. \emph{IEEE Trans. Smart Grid}, 2016. \newline {\tt DOI: 10.1109/TSG.2016.2580945}
    
  \bibitem{zhou}
  D. Zhou, Q. Hu, C. J. Tomlin. Model comparison of a data-driven
and a physical model for simulating HVAC systems. \emph{Amer. Contr. Conf.}, Seattle, 2017.

	
	\bibitem{csm}
	M. Fliess, C. Join. Model-free control.
	\emph{Int. J. Contr.}, 86, 2228-2252, 2013.

\bibitem{pde}
S.D. Howell, P. V. Johnson, P. W. Duck. A rapid PDE-based optimization methodology for temperature control and other mixed stochastic and deterministic systems. \emph{Energ. Buildings}, 43, 1523-1530, 2011.
	
\bibitem{astrom}
	K. J. {\AA}str\"om, T. H\"agglund. \emph{Advanced PID
		Control}, Instrument Soc. Amer., 2006.
	

\bibitem{alinea}
H. Aboua\"{\i}ssa, M. Fliess, C. Join.
On ramp metering: Towards a better understanding of ALINEA via model-free control.
\emph{Int. J. Contr.}, 90, 1018-1026, 2017. 

\bibitem{eccDenmark}
O. Bara, M. Fliess, C. Join, J. Day, S. M. Djouadi. Model-free immune therapy: A control approach to acute inflammation. 
\newblock \emph{Europ. Contr. Conf.}, Aalborg, 2016.	

\bibitem{buc}
H. Wang, X. Ye, Y. Tian, G. Zheng, N. Christov. Model-free based terminal SMC of quadrotor attitude and position. \emph{IEEE Trans. Aerosp. Electron. Syst.}, 52, 2519-2528, 2016.

\bibitem{portugal}
J. Pardelhas, M. Silva, L. Mendon\c{c}a, L. Baptista. Speed control of an experimental pneumatic engine. P. Garrido, F. Soares, A.P. Moreira (Eds): \emph{Controlo 2016}, Lect. Notes Electr. Engin. 402, Springer, 2016.

\bibitem{maa}
S. Maalej, A. Kruszewski, L. Belkoura. Robust control for continuous LPV system with restricted-model-based control. \emph{Circ. Syst. Sign. Process.}, 36, 2499-2520, 2017.


	\bibitem{panc}
T. MohammadRidha, M. Ait-Ahmed, L. Chailloux, M. Krempf, I. Guilhem, J. Y. Poirier, C. H. Moog. Model free iPID control for glycemia regulation of type-1 diabetes. \emph{IEEE Trans. Biomed. Engin.}, 2017. \newline {\tt DOI: 10.1109/TBME.2017.2698036}
	
	\bibitem{its}
	L. Menhour, B. d'Andr\'{e}a-Novel, M. Fliess, D. Gruyer, H. Mounier.  An efficient model-free setting for longitudinal and lateral vehicle control. Validation through the interconnected pro-SiVIC/RTMaps prototyping platform. 
\emph{IEEE Trans. Intel. Transport. Syst.}, 2017. \newline  {\tt DOI: 10.1109/TITS.2017.2699283}

\bibitem{kr}
J. -H. Jeong, D.-H. Lee, Mi. Kim, W.-H. Park, G.-S. Byun. The study of the electromagnetic robot with a four-wheel drive and applied I-PID system. \emph{J. Electr. Eng. Technol.}, 12, 2017.

\bibitem{toulon}
	F. Lafont, J. -F. Balmat, N. Pessel, M. Fliess.
	\newblock A model-free control strategy for an experimental greenhouse with an application to fault accommodation.
	\newblock \emph{Comput. Electron. Agricult.}, 110, 139-149, 2015.
	
	\bibitem{bldg}
H. Aboua\"{\i}ssa, O. Alhaj Hasan, C. Join, M. Fliess, D. Defer. Energy saving for building heating
via a simple and efficient model-free control design: First steps with computer simulations. \emph{Submitted}.
	


	\bibitem{hao2014}
	H. Hao, Y. Lin, A.-S.  Kowli, P. Barooah, S. Meyn.
	\newblock Ancillary service to the grid through control of fans in commercial building HVAC systems.
	\newblock \emph{IEEE Trans. Smart Grid}, 5, 2066-2074, 2014
	
	\bibitem{balanda2014}
	M. Balandat, F. Oldewurtel, M. Chen, C. Tomlin.
	\newblock Contract design for frequency regulation by aggregations of commercial buildings.
	\newblock \emph{52nd  Allerton Conf. Commun. Contr. Comput.}, Monticello, 2014.
	
	\bibitem{siaap}
C. Join, J. Bernier, S. Mottelet, M. Fliess, S. Rechdaoui-Gu\'erin, S. Azimi, V. Rocher. A simple and efficient feedback control strategy for wastewater denitrification. \emph{20th World IFAC Congr.}, Toulouse, 2017.

\bibitem{tn}
C. D. Nelson. \emph{Optimal Control of Energy Efficient Buildings}. Master Thesis, Univ. Tennessee, Knoxville, 2014.

\bibitem{zheng}
L. Zheng, L. Cai. A distributed demand response control strategy using Lyapunov optimization. \emph{IEEE Trans. Smart Grid}, 5, 2075-2083, 2014.

\bibitem{swarm}
L. A. Hurtado, P. H. Nguyen, W.L. Kling. Smart grid and smart building inter-operation using agent-based particle swarm optimization. \emph{Sustain. Energ. Grids Netw.}, 2, 32-40, 2015.

\bibitem{sira1}
M. Fliess, H. Sira-Ram\'{\i}rez.
\newblock An algebraic framework for linear identification.
\newblock \emph{ESAIM Contr. Optimiz. Calc. Variat.},
9, 151-168, 2003.

\bibitem{sira2}
M. Fliess, H. Sira-Ram\'{\i}rez.
\newblock Closed-loop parametric identification for
continuous-time linear systems via new algebraic techniques.
H. Garnier \& L. Wang (Eds): \emph{Identification of
	Continuous-time Models from Sampled Data}, Springer,
pp. 362-391, 2008.

\bibitem{sira3}
H. Sira-Ram\'{\i}rez, C. Garc\'{\i}a-Rodr\'{\i}guez, J. Cort\`{e}s-Romero, A. Luviano-Ju\'{a}rez.  \emph{Algebraic Identification and Estimation Methods in Feedback Control Systems}. Wiley, 2014.

\bibitem{yosida}
K. Yosida. \emph{Operational Calculus} (translated from the
Japanese). Springer, 1984.

	\bibitem{gwerder2005}
M. Gwerder, J. T{\"o}dtli.
\newblock Predictive control for integrated room automation.
\newblock 	\emph{8th REHVA World Congr. Bldg Techno. -- CLIMA 2005}, Lausanne, 2005.

\bibitem{oldewurtel2010}
F. Oldewurtel, A. Parisio,  C. N. Jones, M. Morari, D. Gyalistras, M. Gwerder, V. Stauch, B. Lehmann, K. Wirth.
\newblock Energy efficient building climate control using stochastic model predictive control and weather predictions.
\newblock \emph{Amer. Contr. Conf.}, Baltimore, 2010.

\bibitem{aut}
R. Gondhalekar, F. Oldewurtel, C. N. Jones. Least-restrictive robust periodic model predictive control applied to room temperature regulation. \emph{Automatica}, 49, 2760-2766, 2013.

		
\bibitem{nice}
C. Join, F. Chaxel, M. Fliess. ``Intelligent'' controllers on cheap and small programmable devices.
\emph{2nd Int. Conf. Contr. Fault-Tolerant Syst.}, Nice, 2013. 
		
	
\end{thebibliography}
\end{document}